\documentclass{birkjour}
\usepackage{amsmath,amssymb,amsthm,enumerate}
\usepackage[english]{babel}
\usepackage[T2A]{fontenc}
\usepackage[all]{xy}
 \newtheorem{thm}{Theorem}[section]

 \theoremstyle{definition}
 
 \theoremstyle{remark}

 \numberwithin{equation}{section}

\begin{document}

%
%
%
%
%
%
%
%
%

\title[On the Asymptotic Behavior of Ring Q-homeomorphisms with Respect to P-modulus]
 {On the Asymptotic Behavior \\ of Ring $Q$-homeomorphisms \\ with Respect to $P$-modulus}

\author[Ruslan Salimov]{Ruslan Salimov}

\address{%
01024, Ukraine\\
Kiev-4,\\
3, Tereschenkivska st. }
\email{ruslan.salimov1@gmail.com}

\author{Bogdan Klishchuk}
\address{%
01024, Ukraine\\
Kiev-4,\\
3, Tereschenkivska st. } \email{kban1988@gmail.com}


\date{January 1, 2004}

\begin{abstract}
We study the behavior at infinity of ring $Q$-homeomorphisms with
respect to $p$-modulus for $p>n$.
\end{abstract}

\maketitle
\section{Introduction} Let us recall some definitions, see \cite{Vai}.
Let $\Gamma$ be a family of curves $\gamma$ in $\Bbb R^{n}$,
$n\geqslant2$. A Borel measurable function $\rho:{\Bbb
R^{n}}\to[0,\infty]$ is called {\it admissible} for $\Gamma$, (abbr.
$\rho\in{\rm adm}\,\Gamma$), if
$$
 \int\limits_{\gamma}\rho(x)\,ds\
\geqslant\ 1
$$
for any curve $ \gamma\in\Gamma$. Let $p\in (1,\infty)$.

The quantity
$$
M_p(\Gamma)\ =\ \inf_{\rho\in\mathrm{adm}\,\Gamma}\int\limits_{{\Bbb
R^{n}}}\rho^p(x)\,dm(x)\,
$$
is called {\it $p$--modulus} of the family $\Gamma$.

For arbitrary sets $E$, $F$ and $G$ of $\Bbb R^{n}$ we denote by
$\Delta(E, F, G)$  a set of all continuous curves $\gamma: [a, b]
\rightarrow \Bbb R^{n}$ that connect $E$ and $F$ in $G$, i.e., such
that $\gamma(a) \in E$, $\gamma(b) \in F$ and $\gamma(t) \in G$ for
$a < t < b$.

Let $D$ be a domain in $\Bbb R^{n}$, $n\geqslant2$, $x_0\in D$ and
$d_0 = {\rm dist}(x_{0},
\partial D)$. Set

$$
\mathbb{A}(x_{0}, r_{1}, r_{2}) = \{ x \in \Bbb R^{n}: r_{1} < |x-
x_{0}| < r_{2}\} \,,
$$

$$
S_{i}=S(x_0,r_{i})=\{x\in \Bbb R^{n}: \,  |x-x_0| = r_{i}\}\,, \quad
i = 1,\,2\,.
$$

Let a function $Q: D\rightarrow [0,\infty]$ be Lebesgue measurable.
We say that a homeomorphism $f: D \rightarrow \Bbb R^{n}$ is ring
$Q$-homeomorphism with respect to $p$-modulus at $x_0 \in D$ if the
relation

$$
M_p(\Delta(fS_{1}, fS_{2}, fD))\ \leqslant\
\int\limits_{{\mathbb{A}}}Q(x)\,\eta^{p}(|x-x_{0}|)\,dm(x)\,
$$
holds for any ring $\mathbb{A} = \mathbb{A}(x_{0}, r_{1}, r_{2})$\,,
$0< r_{1} < r_{2} < d_0$, $d_0 = {\rm dist}(x_{0},
\partial D)$, and for any measurable function
$\eta: (r_{1}, r_{2}) \rightarrow [0,\infty]$ such that

$$
\int\limits_{r_{1}}^{r_{2}} \eta(r)\, dr = 1\,.
$$
\medskip

The theory of $Q$-homeomorphisms for $p=n$ was studied in works
\cite{RS}--\cite{S1}, for $1<p<n$ in works \cite{G1}--\cite{S6} and
for $p>n$ in works \cite{SK1}--\cite{SK4}, see also \cite{CR1},
\cite{CR2}.

Denote by $\omega_{n-1}$ the area of the unit sphere
$\mathbb{S}^{n-1} = \{x\in \Bbb R^{n}: \, |x| = 1\}$ in $\Bbb R^{n}$
and by $q_{x_0}(r)= \frac{1}{\omega_{n-1}\,
r^{n-1}}\int\limits_{S(x_0, r)}Q(x)\, d\mathcal{A}$ the integral
mean over the sphere $S(x_0,r)=\{x\in \Bbb R^{n}: \, |x-x_0| =
r\}$\,, here $d\mathcal{A}$ is the element of the surface area.

\medskip

Now we formulate a criterion which guarantees for a homeomorphism to
be the ring $Q$-homeomorphisms with respect to $p$-modulus for $p>1$
in $\Bbb R^{n}$, $n\geqslant2$.

\medskip

{\bf Proposition~1.} {\it Let $D$ be a domain in $\Bbb R^{n}$,
$n\geqslant 2$, and let $Q:D\to[0,\infty]$ be a Lebesgue measurable
function such that $q_{x_0}(r) \neq \infty$ for a.e. $r\in(0,
d_0)$,\, $d_0 = {\rm dist}(x_{0}, \partial D)$. A homeomorphism $f:
D \rightarrow \Bbb R^{n}$ is ring $Q$-homeomorphism with respect to
$p$-modulus at a point $x_{0} \in D$ if and only if the quantity

$$
M_{p}\left(\Delta(fS_{1}, fS_{2}, f\mathbb{A})\right) \leqslant
\frac{\omega_{n-1}}{\left(\int\limits_{r_{1}}^{r_{2}}\frac{dr}{r^\frac{n-1}{p-1}\,q_{x_0}^{\frac{1}{p-1}}(r)}\right)^{p-1}}\,
$$
holds for any $0 < r_{1} < r_{2} < d_0$ } (see \cite{S2}, Theorem
2.3).

\medskip

Following the paper \cite{MRV}, a pair $\mathcal{E}=(A,C)$ where
$A\subset\Bbb R^{n}$ is an open set and $C$ is a nonempty compact
set contained in $A$, is called {\it condenser}. We say that a
condenser $\mathcal{E}=(A,C)$ lies in a domain $D$ if $A\subset D$.
Clearly, if $f:D\to\Bbb R^{n}$ is a homeomorphism and
$\mathcal{E}=(A,C)$ is a condenser in $D$ then $(fA,fC)$ is also
condenser in $fD$. Further, we denote $f\mathcal{E}=(fA,fC)$.

\medskip

Let $\mathcal{E}=(A,C)$ be a condenser. Denote by $\mathcal{C}_0(A)$
a set of continuous functions $u:A\to\mathbb{R}^1$ with compact
support. Let $\mathcal{W}_0(\mathcal{E})=\mathcal{W}_0(A,C)$ be a
family of nonnegative functions $u:A\to\mathbb{R}^1$ such that 1)
$u\in \mathcal{C}_0(A)$, 2) $u(x)\geqslant1$ for $x\in C$ and 3) $u$
belongs to the class ${\rm ACL}$ and

$$
\vert\nabla u\vert= \left(\sum \limits_{i=1}^{n}\left(\frac{\partial
u}{\partial x_i} \right)^{2}\right)^{\frac{1}{2}}\,.
$$
For $p \geqslant1$ the quantity
$$
{\rm cap}_p\,\mathcal{E}={\rm cap}_p\,(A,C)=\inf\limits_{u\in
\mathcal{W}_0(\mathcal{E})}\, \int\limits_{A}\,\vert\nabla
u\vert^p\,dm(x)\,
$$
is called {\it $p$-capacity} of the condenser $\mathcal{E}$. It is
known that for $p>1$
\begin{equation}\label{EMC}
{\rm cap}_p\,\mathcal{E}=M_p(\Delta(\partial A,\partial C;
A\setminus C)),
\end{equation}
see in (\cite{Sh}\,,Theorem~1). For $p>n$ the inequality
\begin{equation}\label{eqks2.8} {\rm cap}_{p}\,  (A,C) \geqslant n\, \Omega_n^{\frac{p}{n}}\, \left(\frac{p-n}{p-1}\right)^{p-1}
\, \left[m^{\frac{p-n}{n(p-1)}}(A) -
m^{\frac{p-n}{n(p-1)}}(C)\right]^{1-p}\,
\end{equation}
holds where $\Omega_n$ is a volume of the unit ball in
$\mathbb{R}^{n}$ (see, e.g., the inequality 8.7 in \cite{Maz}).

\medskip

\section{Main results} Now we consider the main result of our paper
on the behavior at infinity of ring $Q$-homeomorphisms with respect
to $p$-modulus for $p>n$. The case $p=n$ was studied in the work
\cite{SS1}. Let

$$
 L(x_0, f, R) = \sup\limits_{|x - x_{0}|\leqslant R}\, |f(x) - f(x_{0})|\,.
$$

\medskip

\begin{thm}[Main Theorem]
Suppose that $f:\Bbb R^{n} \rightarrow \Bbb R^{n}$ is a ring
$Q$-homeomorphism with respect to $p$-modulus at a point $x_0$ with
$p>n$ where $x_0$ is some point in $\Bbb R^{n}$ and for some numbers
$r_0 >0$, $K>0$ the condition

\begin{equation}
\label{b2H} q_{x_{0}}(t) \leqslant K\,t^{\alpha}\,
\end{equation}
holds for a.e. $t\in [r_0, +\infty)$. If $\alpha\in [0, p-n)$ then

$$
 \varliminf\limits_{R \rightarrow \infty}\, \frac{L(x_0, f, R)}{R^{\frac{p-n-\alpha}{p-n}}} \geqslant K^{\frac{1}{n-p}}\,
  \left(\frac{p-n}{p-n-\alpha}\right)^{\frac{p-1}{p-n}} > 0\,.
$$
If $\alpha = p-n$ then

$$
 \varliminf\limits_{R \rightarrow \infty}\, \frac{L(x_0, f, R)}{\left(\ln R\right)^{\frac{p-1}{p-n}}} \geqslant K^{\frac{1}{n-p}}\,
  \left(\frac{p-n}{p-1}\right)^{\frac{p-1}{p-n}} > 0\,.
$$
\end{thm}
\medskip

\begin{proof}
Consider a condenser $\mathcal{E}=\left(A, C\right)$ in $\Bbb
R^{n}$, where $A=\{x\in \Bbb R^{n}: |x-x_0|< R\}$, $C=\{x\in \Bbb
R^{n}: |x-x_0|\leqslant r_{0}\}$, $0< R< r_{0} <\infty$. Then
$f\mathcal{E} = \left(fA,fC\right)$ is a ringlike condenser in $\Bbb
R^{n}$ and by (\ref{EMC}) we have equality

$$
{\rm cap}_{p}\,  f\mathcal{E}=
\mathrm{M}_{p}\left(\Delta(\partial fA, \partial fC; f(A\setminus
C))\right).
$$
Due to the inequality (\ref{eqks2.8})
$$
{\rm cap}_{p}\,  (fA,fC) \geqslant n\, \Omega_n^{\frac{p}{n}}\, \left(\frac{p-n}{p-1}\right)^{p-1}
\, \left[m^{\frac{p-n}{n(p-1)}}(fA) -
m^{\frac{p-n}{n(p-1)}}(fC)\right]^{1-p}\,
$$
we obtain
\begin{equation}
\label{d3} {\rm cap}_{p}\,  (fA,fC) \geqslant n\,
\Omega_{n}^{\frac{p}{n}}\, \left(\frac{p-n}{p-1}\right)^{p-1}\,
\left[m(fA) \right]^{\frac{n-p}{n}}\,.
\end{equation}
On the other hand, by Proposition~1, one gets

\begin{equation}
\label{d4} {\rm cap}_{p}\,  (fA,fC) \leqslant
\frac{\omega_{n-1}}{\left(\int\limits_{r_0}^{R}\frac{dt}{t^\frac{n-1}{p-1}\,q_{x_0}^{\frac{1}{p-1}}(t)}
\right)^{p-1}}\,.
\end{equation}
Combining the inequalities (\ref{d3}) and (\ref{d4}), we obtain

$$
n\, \Omega_{n}^{\frac{p}{n}}\, \left(\frac{p-n}{p-1}\right)^{p-1}\,
\left[m(fA) \right]^{\frac{n-p}{n}} \leqslant
\frac{\omega_{n-1}}{\left(\int\limits_{r_0}^{R}\frac{dt}{t^\frac{n-1}{p-1}\,q_{x_0}^{\frac{1}{p-1}}(t)}
\right)^{p-1}}\,.
$$
Due to $\omega_{n-1} = n\, \Omega_{n}$, the last inequality can be
rewritten as

\begin{equation}\label{a7}
\Omega_{n}^{\frac{p}{n}- 1}\, \left(\frac{p-n}{p-1}\right)^{p-1}\,
\left[m(fA) \right]^{\frac{n-p}{n}} \leqslant
\left(\int\limits_{r_0}^{R}\frac{dt}{t^\frac{n-1}{p-1}\,q_{x_0}^{\frac{1}{p-1}}(t)}
\right)^{1-p}\,.
\end{equation}

Consider a case when $\alpha \in [0, p-n)$. Then from the condition
(\ref{b2H}) the estimate

$$
\Omega_{n}^{\frac{p}{n}- 1}\, \left(\frac{p-n}{p-1}\right)^{p-1}\,
\left[m(fA) \right]^{\frac{n-p}{n}} \leqslant K\,
\left(\frac{p-n-\alpha}{p-1}\right)^{p-1}\,
\left(R^{\frac{p-n-\alpha}{p-1}} -
r_0^{\frac{p-n-\alpha}{p-1}}\right)^{1-p}\,
$$
holds. Therefore

\begin{equation}\label{a9}
m(fB(x_0, R)) \geqslant \Omega_{n}\, K^{\frac{n}{n-p}}\,
\left(\frac{p-n}{p-n-\alpha}\right)^{\frac{n(p-1)}{p-n}}\,
\left(R^{\frac{p-n-\alpha}{p-1}} -
r_0^{\frac{p-n-\alpha}{p-1}}\right)^{\frac{n(p-1)}{p-n}}\,.
\end{equation}
Due to
\begin{equation}\label{a10}
m(fB(x_0, R)) \leqslant \Omega_{n}\, L^{n}(x_0, f, R)\,,
\end{equation}
from the inequality (\ref{a9}) we have

$$
L(x_0, f, R) \geqslant K^{\frac{1}{n-p}}\,
\left(\frac{p-n}{p-n-\alpha}\right)^{\frac{p-1}{p-n}}\,
\left(R^{\frac{p-n-\alpha}{p-1}} -
r_0^{\frac{p-n-\alpha}{p-1}}\right)^{\frac{p-1}{p-n}}\,.
$$
Dividing the last inequality by $R^{\frac{p-n-\alpha}{p-n}}$ and
taking the lower limit for $R \rightarrow\infty$, we conclude

$$
 \varliminf\limits_{R \rightarrow \infty}\, \frac{L(x_0, f, R)}{R^{\frac{p-n-\alpha}{p-n}}} \geqslant K^{\frac{1}{n-p}}\,
  \left(\frac{p-n}{p-n-\alpha}\right)^{\frac{p-1}{p-n}}\,.
$$

Now we consider a case when $\alpha= p-n$. Then from (\ref{a7}) we
get

$$
\Omega_{n}^{\frac{p}{n}- 1}\, \left(\frac{p-n}{p-1}\right)^{p-1}\,
\left[m(fA) \right]^{\frac{n-p}{n}} \leqslant K\,
\left(\ln\frac{R}{r_0}\right)^{1-p}\,.
$$
Therefore

$$
m(fB(x_0, R)) \geqslant \Omega_{n}\, K^{\frac{n}{n-p}}\,
\left(\frac{p-n}{p-1}\right)^{\frac{n(p-1)}{p-n}}\,
\left(\ln\frac{R}{r_0}\right)^{\frac{n(p-1)}{p-n}}\,.
$$
Due to the estimate (\ref{a10}) we obtain

$$
L(x_0, f, R) \geqslant K^{\frac{1}{n-p}}\,
\left(\frac{p-n}{p-1}\right)^{\frac{p-1}{p-n}}\,
\left(\ln\frac{R}{r_0}\right)^{\frac{p-1}{p-n}}\,.
$$
Finally, dividing the last inequality by $\left(\ln
R\right)^{\frac{p-1}{p-n}}$ and taking the lower limit for $R
\rightarrow\infty$, we conclude

$$
 \varliminf\limits_{R \rightarrow \infty}\, \frac{L(x_0, f, R)}{\left(\ln R\right)^{\frac{p-1}{p-n}}} \geqslant K^{\frac{1}{n-p}}\,
  \left(\frac{p-n}{p-1}\right)^{\frac{p-1}{p-n}}\,.
$$
This completes the proof of Main Theorem.
\end{proof}

\bigskip
Let us consider some examples.

{\it Example 2.1.} Let  $f_1:\mathbb{R}^{n} \to \mathbb{R}^{n}$,
where

$$
f_1(x)=\begin{cases} K^{\frac{1}{n-p}} \left(\frac{p-n}{p-n-\alpha}\right)^{\frac{p-1}{p-n}}\,  |x|^{\frac{p-n-\alpha}{p-n}}\, \frac{x}{|x|} \, ,&  x \neq 0\\
0,&  x=0 \,.\end{cases}\,
$$
It can be easily seen that $\lim\limits_{x \rightarrow \infty}\,
\frac{|f(x)|}{|x|^{\frac{p-n-\alpha}{p-n}}} = K^{\frac{1}{n-p}}\,
\left(\frac{p-n}{p-n-\alpha}\right)^{\frac{p-1}{p-n}}$. Let us show
that the mapping $f_1$ is a ring $Q$-homeomorphism with respect to
$p$-modulus with the function $Q(x) = K\,|x|^{\alpha}$ at the point
$x_{0}= 0$. Clearly, $q_{x_{0}}(t) = K\,t^{\alpha}$. Consider a ring
$\mathbb{A}(0, r_{1}, r_{2})$, $0< r_{1} < r_{2} < \infty$. Note
that the mapping $f_{1}$ maps the ring $\mathbb{A}(0, r_{1}, r_{2})$
onto the ring $\widetilde{\mathbb{A}}(0, \widetilde{r}_{1},
\widetilde{r}_{2})$, where
$$
\widetilde{r}_{i} = K^{\frac{1}{n-p}}\,
\left(\frac{p-n}{p-n-\alpha}\right)^{\frac{p-1}{p-n}}\,
r_{i}^{\frac{p-n-\alpha}{p-n}}\,, \quad i = 1,\,2.
$$
Denote by $\Gamma$ a set of all curves that join the spheres
$S(0,r_{1})$ and $S(0,r_{2})$ in the ring $\mathbb{A}(0, r_{1},
r_{2})$. Then one can calculate $p$-modulus of the family of curves
$f_{1}\Gamma$ in implicit form:

$$
\mathrm{M}_{p}(f_{1}\Gamma) = \omega_{n-1}\,
\left(\frac{p-n}{p-1}\right)^{p-1}\,
\left(\widetilde{r}_{2}^{\frac{p-n}{p-1}} -
\widetilde{r}_{1}^{\frac{p-n}{p-1}}\right)^{1-p}\,
$$
(see, e.g., the relation (2) in \cite{Ge}). Substituting in the
above equality the values $\widetilde{r}_{1}$ and
$\widetilde{r}_{2}$, defined above, one gets
 $$
 \mathrm{M}_{p}(f_{1}\Gamma) =
 \omega_{n-1}\, K\, \left(\frac{p-n-\alpha}{p-1}\right)^{p-1}\, \left(r_{2}^{\frac{p-n-\alpha}{p-1}} -
 r_{1}^{\frac{p-n-\alpha}{p-1}}\right)^{1-p}\,.
$$
Note that the last equality can be written by
$$
\mathrm{M}_{p}(f_{1}\Gamma) =
\frac{\omega_{n-1}}{\left(\int\limits_{r_{1}}^{r_{2}}\frac{dt}{t^\frac{n-1}{p-1}\,q_{x_0}^{\frac{1}{p-1}}(t)}\right)^{p-1}}\,,
$$
where $q_{x_0}(t) = K\,t^{\alpha}$.

Hence, by Proposition 1, the homeomorphism $f_{1}$ is a ring
$Q$-homeomorphism with respect to $p$-modulus for $p>n$ with the
function $Q(x) = K\,|x|^{\alpha}$ at the point $x_{0}= 0$.

\medskip

{\it Example 2.2.} Let $\alpha= p-n$ and $f_2:\mathbb{R}^{n} \to
\mathbb{R}^{n}$, where

$$
f_2(x)=\begin{cases} K^{\frac{1}{n-p}} \left(\frac{p-n}{p-1}\right)^{\frac{p-1}{p-n}}\,  \left(\ln|x|\right)^{\frac{p-1}{p-n}}\, \frac{x}{|x|} \, ,&  x \neq 0\\
0,&  x=0 \,.\end{cases}\,
$$
It can be easily seen that $\lim\limits_{x \rightarrow \infty}\,
\frac{|f(x)|}{\left(\ln|x|\right)^{\frac{p-1}{p-n}}} =
K^{\frac{1}{n-p}}\, \left(\frac{p-n}{p-1}\right)^{\frac{p-1}{p-n}}$.
By analogy with Example 2.1, we can show that the mapping $f_2$ is a
ring $Q$-homeomorphism with respect to $p$-modulus with the function
$Q(x) = K\,|x|^{p-n}$.

{\it Remark 2.1.} Examples 2.1 and 2.2 show that the estimates in
Main Theorem are exact, i.e. are attained on the above mappings.

This work was supported by the budget program ``Support of the
develo\-pment of priority trends of scientific
researches''\,(KPKVK\,6541230).

\end{document}